\documentclass[11pt]{article}
\usepackage{amsfonts}
\usepackage{mathrsfs,color,amscd,amssymb,enumerate,amsthm,amsmath,bm,graphicx,psfrag,subfigure}
\usepackage{latexsym}
\usepackage{float,fancybox,shapepar,setspace,hyperref}
\usepackage{pgf,tikz}

\makeatletter
\def\leftharpoonfill@{\arrowfill@\leftharpoonup\relbar\relbar}
\def\rightharpoonfill@{\arrowfill@\relbar\relbar\rightharpoonup}
\newcommand\rbjt{\mathpalette{\overarrow@\rightharpoonfill@}}
\newcommand\lbjt{\mathpalette{\overarrow@\leftharpoonfill@}}
\makeatother

\makeatletter

\renewcommand{\@seccntformat}[1]{{\csname the#1\endcsname}{\normalsize .}\hspace{.5em}}
\makeatother

\def \[{\begin{equation}}
\def \]{\end{equation}}

\def \ss {\subseteq}

\def\bqed{ \hfill $\blacksquare$}

\newtheorem{thm}{Theorem}

\newtheorem{claim}{Claim}

\newtheorem{lem}{Lemma}

\newtheorem{pb}{Problem}

\newtheorem{conj}{Conjecture}

\newtheorem{definition}{Definition}

\usetikzlibrary{arrows}
\begin{document}

\title{Oriented  diameter of graphs with given domination number}
\author{Xiaolin WANG$^1$, Yaojun CHEN$^{2,}$\footnote{Corresponding author. Email: yaojunc@nju.edu.cn}\\
\small{$^1$School  of  Mathematics  and  Statistics,  Fuzhou  University,  Fuzhou  350108,  P.R. CHINA}\\
 \small{$^2$School of Mathematics, Nanjing University, Nanjing
210093, P.R. CHINA}}
\date{}
\maketitle
\begin{abstract}

Let $G$ be a connected bridgeless graph with domination number $\gamma$. The oriented diameter  (strong  diameter) of $G$ is the smallest integer $d$ for which $G$ admits a strong orientation with diameter (strong diameter) $d$. Kurz and L\"atsch (2012) conjectured the oriented diameter of $G$ is at most $\lceil \frac{7\gamma+1}{2}\rceil$ and the bound is sharp. 
In this paper, we confirm the conjecture by induction on $\gamma$ through contracting an unavoidable alternative subgraph, which holds potential for future applications. Moreover, we show the  oriented strong  diameter of $G$ is at most $7\gamma-1$ by using the same recursive structure, and the bound is best possible.

\vskip 2mm

\end{abstract}
{\setcounter{section}{0}

\section{Introduction}\setcounter{equation}{0}
\vskip 2mm
Let $G=(V(G),E(G))$ be a  simple undirected connected graph.
For any $u,v\in V(G)$, the distance $d(u,v)$ is the length of a shortest $(u,v)$-path connecting $u$ and $v$, and  the  \textit{diameter} of $G$ is defined as $diam(G)=max\{d(u,v)\,|\,u,v \in V(G)\}$. The degree of a vertex $v$ in $G$, denoted by $d_G(v)$, is the number of vertices adjacent to $v$. A \textit{bridge} of  a graph $G$ is an edge whose removal disconnects $G$. A graph $G$ is {\it  bridgeless} if it has no bridge. A maximal 2-connected subgraph of a graph is referred as a \emph{block}. For a connected subgraph $H$ of a graph $G$, $G/H$ denotes the resulting graph by contracting $H$ into a vertex. 
 An orientation $\rbjt{G}$ of a graph $G$ is a digraph obtained from it by assigning a direction to each edge. 
 An orientation $\rbjt{G}$ is \emph{strong} if there is a directed path from $u$ to $v$ for any $u,v\in V(\rbjt{G})$.
The directed distance $\partial(u,v)$ is the length of a shortest directed path from $u$ to $v$ in $\rbjt{G}$. Let $\theta(u,v)=max\{\partial(u,v),\partial(v,u)\}$.
If $\rbjt{G}$ is strong, we define the diameter of $\rbjt{G}$ as $diam(\rbjt{G})=max\{\theta(u,v)\,|\,u,v \in V(\rbjt{G})\}$.

In 1939, Robbins \cite{R} gave the well-known result that an undirected connected graph admits a strong orientation if and only if it is bridgeless.  For a  bridgeless graph $G$, 
define the \textit{oriented diameter} as
$$\rbjt{diam}(G)=min\left\{diam(\rbjt{G})~\left|~\rbjt{G} ~\text{is a strong orientation of}~ G\right\}.\right.$$

The notion of oriented diameter was first introduced by Chv\'atal and Thomassen in \cite{CT}, and they considered the relationship between the diameter and oriented diameter in the same paper.
 Let $f(d)$ be the smallest value for which every bridgeless graph $G$ with diameter $d$ admits a strong orientation $\rbjt{G}$ such that  
$diam(\rbjt{G})\leq f(d)$. Chv\'atal and Thomassen  \cite{CT} established the following bounds for $f(d)$.

\begin{thm}(Chv\'atal and Thomassen  \cite{CT})\label{CT}
$\frac{1}{2}d^2+d\leq f(d)\leq 2d^2+2d.$
\end{thm}
By Theorem \ref{CT}, $f(d)=\Theta(d^2)$. Babu et al. \cite{BB} improved the upper bound  as follows, which is smaller than $2d^2+2d$ when $d\geq 8$.
 
 \begin{thm}(Babu et al. \cite{BB} )
$ f(d)\leq1.37d^2+6.97d.$
 \end{thm}

Although we have known that $f(d)=\Theta(d^2)$, it seems difficult to determine the exact value of $f(d)$, even to determine the coefficient of term $d^2$ in general. Indeed, it is very interesting to determine the exact value of $f(d)$ when $d$ is small.  For $d=2$, Chv\'atal and Thomassen  \cite{CT} proved the following.
\begin{thm}(Chv\'atal and Thomassen  \cite{CT})
$f(2)=6$.

\end{thm}

For $d=3$,
Kwok, Liu and West \cite{KLW} proved that $9\leq f(3)\leq 11$ which greatly improved the bounds for $f(3)$ obtained by Theorem \ref{CT}. Later,  Wang and Chen \cite{WC} determined the exact value of  $f(3)$.

\begin{thm}(Wang and Chen \cite{WC})
$f(3)=9.$
\end{thm}

Up to now, only two exact values of $f(d)$ are determined. Furthermore, there are many other researches about oriented diameter:  obtaining sharp upper bounds for some special classes of graphs (\cite{FM}-\cite{KRS}, \cite{LS}, \cite{WCD}),   and establishing  tight upper bounds in terms of other graph parameters, including maximum degree (\cite{CC}, \cite{DGS}), minimum degree (\cite{BD}, \cite{MS}) and so on.

 A dominating set $D$ of a graph $G$ is a vertex set such that any $v \in V(G) \backslash D$ is adjacent to at least one vertex of $D$. The domination number $\gamma(G)$ of a graph $G$ is the minimum cardinality among all dominating sets, which is a fundamental graph parameter to be well studied in graph theory. 
In 2004, 
Fomin, Matamala, Prisner,
and Rapaport  \cite{FM}
started to investigate the relations between the oriented diameter and domination number, and obtained the following. 

\begin{thm}\label{Fomin}
(Fomin, Matamala, Prisner,
and Rapaport  \cite{FM}). Let $G$ be a bridgeless graph with $\gamma(G)=\gamma$.
Then  $\rbjt{diam}(G) \leq 9\gamma$.
\end{thm}

In 2012,  Kurz and  L\"atsch \cite{KL} greatly improved this bound as follows.

\begin{thm}\label{Kurz}
(Kurz and  L\"atsch \cite{KL}). Let $G$ be a bridgeless graph with $\gamma(G)=\gamma$.
Then $\rbjt{diam}(G) \leq 4\gamma$.
\end{thm}
In the same paper, they posed the following conjecture. 

\begin{conj}\label{conjecture}(Kurz and  L\"atsch \cite{KL}). Let $G$ be a bridgeless graph with $\gamma(G)=\gamma$.
Then $\rbjt{diam}(G) \leq \lceil \frac{7\gamma+1}{2}\rceil$.
\end{conj}

If Conjecture \ref{conjecture} is true, the upper bound is tight as can be seen by the graph described as follows \cite{KL}:  For a path $P=u_1\cdots u_{\gamma}$, replace each edge $u_{i-1}u_i$ by the graph as shown in Figure \ref{tight} (1) or (2), and replace $u_1,u_\gamma$ by a triangle as shown in Figure \ref{tight} (3). Let $G$ be the resulting graph. It is not difficult to check that $\{u_1,\ldots,u_\gamma\}$ is a minimum dominating set of $G$, and $\rbjt{diam}(G)=\lceil \frac{7\gamma+1}{2}\rceil$. 
\vskip 2mm
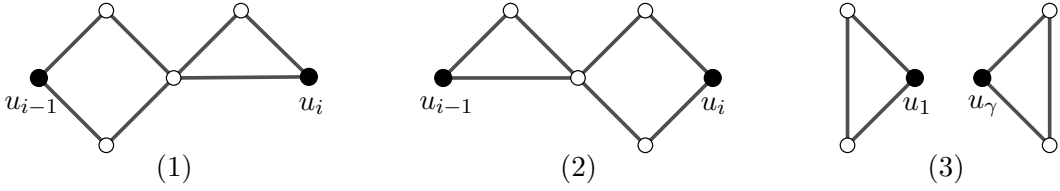
\begin{figure}[!htb]
\centering
\includegraphics[height=0.16\textwidth]{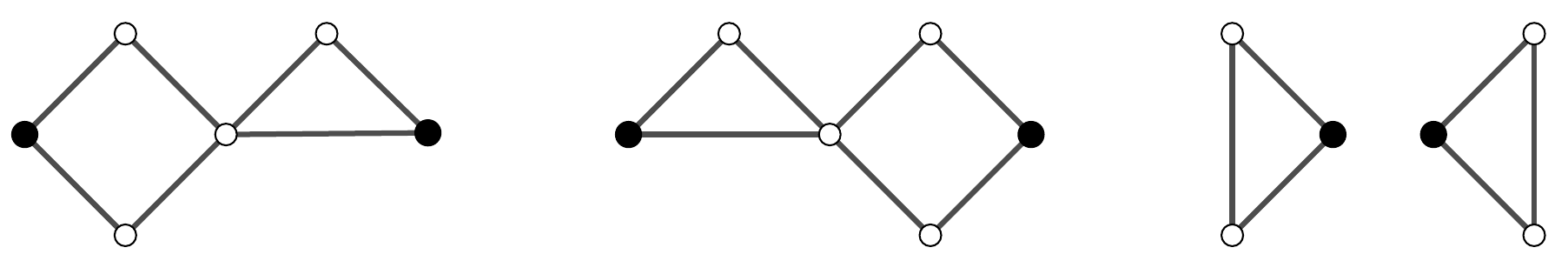}
\put(-341,-4){\makebox(3,3){$(1)$}}
\put(-188,-4){\makebox(3,3){$(2)$}}
\put(-49,-4){\makebox(3,3){$(3)$}}
\put(-395,20){\makebox(3,3){$u_{i-1}$}}
\put(-289,20){\makebox(3,3){$u_{i}$}}
\put(-238,20){\makebox(3,3){$u_{i-1}$}}
\put(-137,20){\makebox(3,3){$u_{i}$}}
\put(-60,20){\makebox(3,3){$u_{1}$}}
\put(-35,20){\makebox(3,3){$u_{\gamma}$}}
\vskip 3mm
\caption{The construction of the extremal graph}
\label{tight}
\end{figure}

In this paper, we confirm that Conjecture \ref{conjecture} is true and obtain the following.

\begin{thm}\label{main}
Let $G$ be a bridgeless graph and $D$ a dominating set of $G$ with $|D|=\gamma$. Then $\rbjt{diam}(G) \leq \lceil \frac{7\gamma+1}{2}\rceil$, and this bound is sharp.
\end{thm}

The main idea for proving Theorem \ref{main} is by induction on $\gamma$. For this purpose, we try to find an appropriate bridgeless subgraph $H$ of $G$ with $|V(H)\cap D|\ge 2$. By induction hypothesis, $G/H$ admits a strong orientation with diameter as expected. Choose an appropriate strong orientation for $H$. Based on the given strong orientations of $G/H$ and $H$, we can get a strong orientation of $G$ as required. The most difficult part for doing this is how to find the appropriate subgraph $H$ for contraction, which is the main goal we discuss in the next section.

If $\rbjt{G}$ is strong,
 the strong distance $sd(u,v)$ is the minimum number of arcs among all strong subgraphs  of $\rbjt{G}$ containing $u$ and $v$. Obviously,   $sd(u,v) \leq \partial(u,v)+\partial(v,u)$. 
 This notion was first introduced by Chartrand, Erwin, Raines and Zhang \cite{CER} in 1999. It is motivated by the definition of directed distance in digraph, which is  lack of symmetry, and by the definition of distance of  $u$ and $v$ in undirected graphs, which can be viewed as the minimum number of edges among all  connected subgraphs  of $G$ containing $u$ and $v$. 
 If $\rbjt{G}$ is strong, we define  the strong diameter of $\rbjt{G}$ as $sdiam(\rbjt{G})=max\{sd(u,v)\,|\,u,v \in V(\rbjt{G})\}$.  For a bridgeless graph $G$, we define the \textit{oriented strong  diameter} as
$$\rbjt{sdiam}(G)=min\left\{sdiam(\rbjt{G})~\left|~\rbjt{G} ~\text{is a strong orientation of}~ G\right\}.\right.$$
The notion of oriented strong  diameter was  introduced by Lai, Chiang,  Lin and Yu  \cite{LCL}. And there are some researches about oriented strong diameter for some special classes of graphs \cite{CG}-\cite{CGZ}, \cite{GH}-\cite{GH2}, \cite{LCL}-\cite{MG2}.

Using the same method and the subgraph structure as in the proof of Theorem \ref{main}, we establish a tight bound for the oriented strong  diameter as follows.

\begin{thm}\label{main2}
Let $G$ be a bridgeless graph and $D$ a dominating set of $G$ with $|D|=\gamma$. Then  $\rbjt{sdiam}(G) \leq 7\gamma-1$, and this bound is sharp.
\end{thm}

It is easy to check that the extremal graph as shown in Figure \ref{tight} for oriented diameter  is also the one for oriented strong diameter.

\section{Unavoidable structure in standard pair (G,D)}

Let $G$ be a bridgeless graph and $D$ a dominating set of $G$ with $|D|=r$. We call the vertices in $D$ as \emph{dominating} vertices and the vertices in $V(G)\backslash D$ as \emph{dominated} vertices. We say that a vertex $v\in V(G)$ is \emph{associated} with an isolated triangle if there exists two vertices $v_1,v_2 \in V(G)$ such that $d_G(v_1)=d_G(v_2)=2$ and $vv_1v_2v$ is a triangle.

The pair $(G,D)$ is called \emph{standard} if $G$ is bridgeless, $D$ is an independent set, any vertex in $V(G)\backslash D$ is dominated by exactly one vertex in $D$, and each dominating vertex is associated with an isolated triangle.
Actually, if $(G,D)$ is not standard, we can always transform it into a standard pair $(G',D)$ by the following three operations.
 \begin{itemize}
\setlength{\itemsep}{1pt}
\setlength{\parsep}{1pt}
\setlength{\parskip}{1pt}
\item 
If there exist  $u,v \in D$ and $uv \in E(G)$, then subdivide $uv$ twice.

\item If $v \in V \backslash D$ is dominated by $t$ vertices $v_1,\ldots,v_t$ in $D$ with $t \geq 2$, then subdivide $vv_i$ once for all $2\le i \le t$. 

\item If $v \in D$ is not associated with an isolated triangle, then add two new vertices $v',v''$ such that $vv'v''$ is a triangle.

\end{itemize}
It is easy to see that such  operations do not decrease its oriented diameter and oriented strong diameter, that is, $\rbjt{diam}(G)\le \rbjt{diam}(G')$ and $\rbjt{sdiam}(G)\le \rbjt{sdiam}(G')$. Because $D$ is still a dominating set in $G'$, it is sufficient to prove Theorems \ref{main} and \ref{main2} for the pair $(G,D)$ being standard. Throughout the rest of this paper, we always assume that $(G,D)$ is a standard pair.

Since each dominating vertex is associated with an isolated triangle and this triangle must be a directed triangle under any strong orientation of $G$, we can get easily the following lemma.

\begin{lem}\label{-2} Let $(G,D)$ be standard and $\rbjt{G}$ be any strong orientation of $G$. If $G$ is not a triangle, then we have $\theta(u,v) \leq diam(\rbjt{G})-2$ for any $u \in V(G)\backslash D$ and $v \in D$.
\end{lem}

Our main idea for proving Theorems \ref{main} and \ref{main2} is by induction on $r$. For this purpose, we need discuss some unavoidable structure in a given standard pair $(G,D)$.
Before discussing, we first give some definitions.

\begin{definition}
A bridgeless subgraph $H_0$ is called \emph{alternative} if: 

(1) $V(H_0)\cap D=\{x_1,\ldots,x_s\}$ and $s\geq 1$; 

(2) $x_i$ dominates exactly two vertices $x_{i1},x_{i2}$ in $H_0$ for any $i\leq s$;

(3) $V(H_0)=\cup_{i=1}^{s} \{x_i,x_{i1},x_{i2}\}$. 
\end{definition}

\begin{definition} For $H\ss G$ and $V(H)\cap D\not=\emptyset$, $S(G/H)$ denotes a graph obtained from $G$ by 
contracting $H$ into a vertex $h$, and then subdivide all the edges $hv$ in $G/H$ as follows. For each edge $hv\in E(G/H)$, denote $h'v$ as the corresponding edge in $E(G)$ and $h' \in V(H)$.

(1) If $h'$ is a dominating vertex in $G$, then 
do not subdivide $hv$ in $S(G/H)$.

(2) If neither $h'$ or $v$ is a dominating vertex in $G$, then subdivide  $hv$ once in $S(G/H)$.

(3) If $h'$ is  dominated by $v$  in $G$, then subdivide  $hv$ twice in $S(G/H)$.

\end{definition}
The graph $S(G/H)$ is illustrated in Figure \ref{cs}, where the dark and white points represent the dominating and  dominated vertices, respectively. It is easy to check  that $(S(G/H),D\cup \{h\}\backslash V(H))$ is standard (add an isolated triangle associated with $h$).
\vskip -1mm
\begin{figure}[!htb]
\centering
\includegraphics[height=0.25\textwidth]{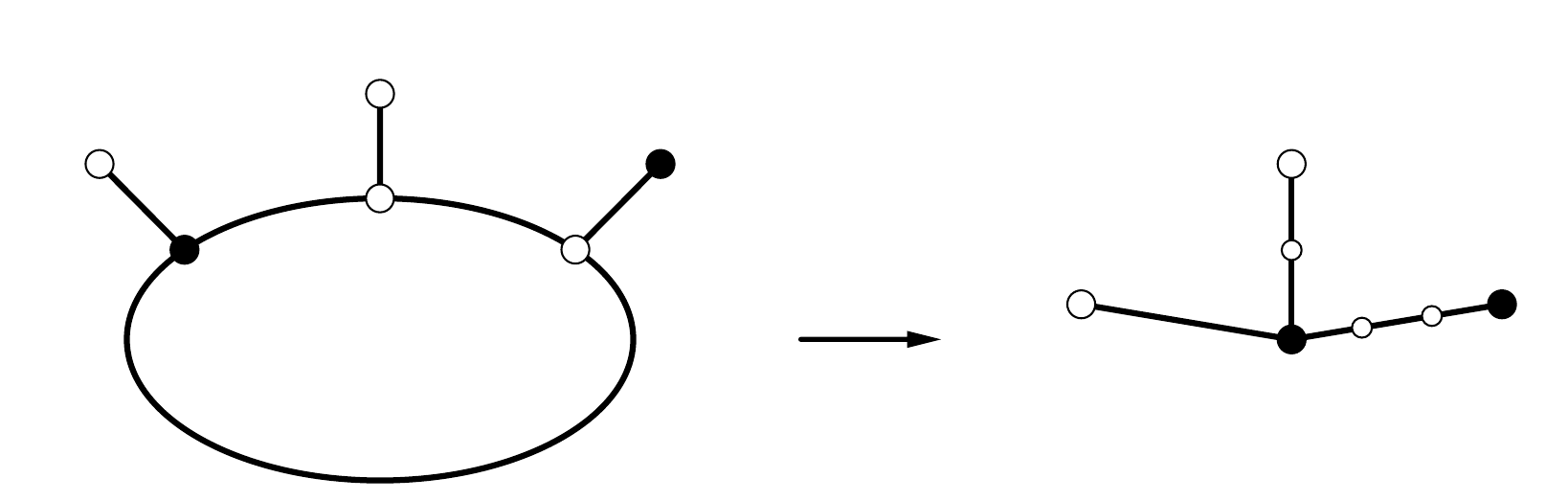}
\put(-58,22){\makebox(3,3){$h$}}
\put(-246,32){\makebox(3,3){$H$}}
\caption{The subgraph $H$ and $S(G/H)$.}
\label{cs}
\end{figure}

\begin{lem}\label{vrh2}  Let $(G,D)$ be a standard pair with $|D|\geq 2$. Then for any $x \in D$, 
  there exists an alternative subgraph $H_0$ containing $x$ and at least two dominating vertices.
\end{lem}

\noindent{\bf{Proof.}}  
In order to find an alternative subgraph $H_0$ as required, we first get a cycle $\mathcal{C}_t$ containing $x$ through a series of cycle contractions recursively  starting from a cycle $C_{v_1}$ in $G=G_1$, and then recover all contracted cycles one by one in reverse order.

Let  $G_1=G$ and $D_1 = D\cap V(G_1)$. Then $(G_1,D_1)$ is  standard. Since $|D_1|\geq 2$,  it is not difficult to see $G_1$ has a maximal path $P_{u_1}$ starting at $x$ and ending at some $u_1 \in D_1$ of type $$x\!\odot\!\odot\!\circledast\cdots \circledast\!\odot\!\odot\!\circledast\cdots \circledast\!\odot\!\odot\! u_1,$$ 
where $\odot$ and $\circledast$ represent a dominated  and a dominating vertex  in $G_1$, respectively.

Now we give the following iteration: suppose that in each step $i\ge 1$,  $(G_i,D_i)$ is  standard, and $P_{u_i}$ is a maximal path starting at $x$ and ending at some $u_{i} \in D_i$ in $G_i$ of type $x\!\odot\!\odot\!\circledast\cdots \circledast\!\odot\!\odot\!\circledast\cdots \circledast\!\odot\!\odot\! u_i.$

Since $G_i$ is bridgeless, there exists a shortest path $Q_{u_i}=u_iq_1\cdots q_\ell v_i$  between $u_i$ and $P_{u_i}-u_i$ other than the last edge of $P_{u_i}$. Because $(G_i,D_i)$ is standard and $P_{u_i}$ is maximal, $Q_{u_i}$ is
of type $u_i\!\odot\! v_i$ or $u_i\!\odot\!\odot v_i$, and $v_i$ is a dominated or dominating vertex, respectively.

Set $C_{v_i}= v_iP_{u_i}u_i Q_{u_i}v_i$. Then $C_{v_i}$ is a cycle in $G_i$. Note that $v_i$ may be  a dominating vertex, or dominated by a vertex in $C_{v_i}$, or dominated by its predecessor in $P_{u_i}$. By the types of $P_{u_i}$ and $Q_{u_i}$, we have the following.

\begin{claim}\label{cvi} $V(C_{v_i})\cap D_i\not=\emptyset$. Moreover, if $v_i$ is dominated by its predecessor in $P_{u_i}$, then every vertex in $V(C_{v_i})\cap D_i$ dominates two vertices in $C_{v_i}$, and for all dominated vertices in $C_{v_i}$, only $v_i$ is dominated by a vertex not in $C_{v_i}$. Otherwise,
$C_{v_i}$ is always alternative.

\end{claim}

If $v_i=x$, then we get an alternative cycle containing two dominating vertices $x$ and $u_i$ in $G_i$, and the iteration stops.
Otherwise, we let $G_{i+1}=S(G_i/C_{v_i})$ and $v_i$ be the contracted vertex in $G_i/C_{v_i}$. It is easy to see that $(G_{i+1},D_{i+1})$ is  standard with $D_{i+1}=(D_i\backslash V(C_{v_i}))\cup\{v_i\}$. Let $P_{i+1}$ be the path in $G_{i+1}$ corresponding to $P_{u_i}$ after contraction and subdivision. Extend $P_{i+1}$
 to a maximal path $P_{u_{i+1}}$  starting at $x$ and ending at some $u_{i+1} \in D_{i+1}$ in $G_{i+1}$ of type $x\!\odot\!\odot\!\circledast\cdots \circledast\!\odot\!\odot\!\circledast\cdots \circledast\!\odot\!\odot\! u_{i+1}.$
Then we go to step $i+1$.

Note that we may add some new vertices after the subdivision operation,  and these vertices have degree 2. Because $v_i\in V(C_{v_i})$ and $d_{G_i}(v_i)\geq 3$, by the construction of $G_{i+1}$, the sum of the degree of vertices with degree at least 3 in $G_{i+1}$ decreases at least 2 comparing to that in $G_i$. Hence, this iteration will stop within finite steps. 
When it ends at step $t$, $v_t=x$ and we obtain cycles $C_{v_{i}}$ in $G_{i}$ for $1\le i \leq t$. By Claim \ref{cvi},  $C_{v_{t}}$ is an alternative subgraph in $G_t$ containing dominating vertices $x$ and $u_t$. 

\vskip 2mm
Set $\mathcal{C}_t=C_{v_{t}}$. We now begin to recover contracted cycles one by one starting at $\mathcal{C}_t$ to get required $H_0$. In this process, the most difficult situation is to recover $C_{v_i}$ in the case when $v_i$ is dominated by its predecessor in $P_{u_i}$. In order to do this, we need to consider a very special path of order four defined as follows.

In the step $j$,  if $v_j$ is dominated by its predecessor $y_j$ in $P_{u_j}$, then after the contraction and subdivision operation,  $v_j$ becomes a dominating vertex, and  $y_jv_j$ becomes $y_jz_{j1}z_{j2}v_j$. Note that $v_j$ may be contracted into $v_k$ for some $k>j $.
For all $k \ge j$, if $y_jz_{j1}z_{j2}v_k$ is a subpath of $P_{u_k}$, then for convenience, let
$y_kz_{k1}z_{k2}v_k=y_jz_{j1}z_{j2}v_k$ and 
refer it as a \emph{bad} $(y_k,v_k)$-path. Obviously, for any such $v_k$, we have only one bad $(y_k,v_k)$-path. Since $z_{j1},z_{j2}$ are two subdivided vertices, by the above iteration, we have the following claim which tells us where the bad path should occur.

\begin{claim}\label{detail} In step $i$, for any bad $(y_j,v_j)$-path $y_jz_{j1}z_{j2}v_j$ with $j<i$, if $v_j \in P_{u_i}$, then $y_jz_{j1}z_{j2}v_j$  must occur in $P_{u_i}$ of the form $P_{u_i}=x\cdots y_jz_{j1}z_{j2}v_j\cdots u_i$.

\end{claim}

We now continue the recover process. Applying Claim \ref{detail} on $i=t$, 
for any bad $(y_j,v_j)$-path $y_jz_{j1}z_{j2}v_j$ with $j<t$, if $v_j \in \mathcal{C}_t$, then $y_jz_{j1}z_{j2}v_j \ss \mathcal{C}_{t}$. Assume that we have got alternative subgraph $\mathcal{C}_{i+1}$ containing $x$ and $|V(\mathcal{C}_{i+1})\cap D_{i+1}|\geq 2$ in $G_{i+1}$ with $i+1\le t$, and $\mathcal{C}_{i+1}$ satisfies the following property: 
 \vskip 2mm
\noindent For any bad $(y_j,v_j)$-path $y_jz_{j1}z_{j2}v_j$ with $j<i+1$,  $y_jz_{j1}z_{j2}v_j \ss \mathcal{C}_{i+1}$ if $v_j \in \mathcal{C}_{i+1}$.\hfill{($*$)}
\vskip 2mm
\noindent By the above iteration, we can see that in the step $i$, $C_{v_{i}}$ is contracted to a dominating vertex $v_i$ in $G_{i+1}$. We first recover all subdivided edges in $G_{i+1}$ and then recover $C_{v_{i}}$ in $G_{i+1}$, and thus result in $G_i$. Assume $\mathcal{C}_{i+1}$ becomes $\mathcal{C}'_i\ss G_i$ after $G_{i+1}$ recovered to $G_i$.
This  $\mathcal{C}_i'$ may be an alternative subgraph. In such a case, $\mathcal{C}_i'$ satisfies ($*$).  Sometimes, this $\mathcal{C}_i'$  is not an alternative subgraph, and we need carefully choose a path from $C_{v_i}$ to delete in $\mathcal{C}_i'$ so as to get an alternative subgraph $\mathcal{C}_i$ satisfying ($*$).
Seeing this, in the current step, $v_i$ is not necessarily lying in $\mathcal{C}_i'$. If $v_i \notin \mathcal{C}_{i}'$, then $\mathcal{C}_i=\mathcal{C}'_i=\mathcal{C}_{i+1}\ss G_i$ is required. Otherwise, 
 we will find a $\mathcal{C}_i$ as required within $\mathcal{C}_i'$ as follows.

Clearly, $v_i$ is of degree 2 in $C_{v_i}$. Since $\mathcal{C}_{i+1}$ is alternative and $v_i$ is a dominating vertex, $v_i$ has degree 2 in $\mathcal{C}_{i+1}$. Thus, $v_i$ has degree 3 or 4  in $\mathcal{C}_i'$, that is, $d(v_i)=4$ in $\mathcal{C}_i'$, or $C_{v_{i}}$ contains two vertices $v_i$ and $v_i'$ which have degree 3 in $\mathcal{C}_i'$. Note that $v_i$ is a dominating vertex in $\mathcal{C}_{i+1}\ss G_{i+1}$, but it may be a dominating or  dominated  vertex in $\mathcal{C}_i'\ss G_i$. If $v_i$ is a dominated vertex in $\mathcal{C}_i'\ss G_i$, then by the construction of $P_{u_i}$ and $C_{v_i}$,  
 $v_i$ is dominated by some vertex in $V(C_{v_i})$, or  by its predecessor $y_i$ in $P_{u_i}$.

 Suppose that $d(v_i)=4$ in $\mathcal{C}_i'$. If $v_i$ is dominated by some vertex in $V(C_{v_i})$, 
then by the subdivision operation, we can check that two neighbors of $v_i$ in $\mathcal{C}_{i+1}$ are two subdivided vertices and are dominated by $v_i$ in $\mathcal{C}_{i+1}$. Moreover, these two vertices are contracted as $\mathcal{C}_{i+1}$ becomes $\mathcal{C}_i'$. 
 When we recover $C_{v_i}$, $v_i$ becomes  a dominated vertex in $G_i$. Thus, combining Claim \ref{cvi}, 
we can check that $\mathcal{C}_i=\mathcal{C}_i'$ is required. 

 If $v_i$ is dominated by its predecessor $y_i$ in $P_{u_i}$,
then by the subdivision operation, $y_iz_{i1}z_{i2}v_i$ is a bad $(y_i,v_i)$-path in $P_{u_{i+1}}$. By ($*$), we can see that 
$y_iz_{i1}z_{i2}v_i \ss \mathcal{C}_{i+1}$. Since $d(v_i)=4$ in $\mathcal{C}_i'$, by the subdivision operation, let $z_{i3}$ be another subdivided vertex in $\mathcal{C}_{i+1}$ dominated by $v_i$. Then, $z_{i1},z_{i2},z_{i3}$ are contracted as $\mathcal{C}_{i+1}$ becomes $\mathcal{C}_i'$. 
When we recover $C_{v_i}$, $v_i$ becomes  a vertex dominated by $y_i\in \mathcal{C}_i'-V(C_{v_i})$. Thus, combining Claim \ref{cvi}, 
we can check that $\mathcal{C}_i=\mathcal{C}_i'$ is required.

  If $v_i \in D_i$, then by the contraction and subdivision operation, $\mathcal{C}_i'=\mathcal{C}_{i+1}\cup C_{v_i}$. Let $\mathcal{C}_i=\mathcal{C}_{i+1}$, then $\mathcal{C}_i$ is alternative containing at least two dominating vertices $x$ and $v_i$. Because $\mathcal{C}_i$ is obtained from $\mathcal{C}_i'$ by deleting some vertices, we need to show that $\mathcal{C}_i=\mathcal{C}_{i+1}$ still satisfies ($*$).   
  If there is a bad $(y_j,v_j)$-path such that $v_j \in \mathcal{C}_i$ 
 and $y_jz_{j1}z_{j2}v_j \nsubseteq \mathcal{C}_{i}$ for some $j<i$, then $v_j$ is contracted to $v_i$ and $y_jz_{j1}z_{j2}v_j \ss C_{v_i}$, which contradicts Claim \ref{detail} because $x$ is in $\mathcal{C}_{i+1}$.  Therefore, $\mathcal{C}_i=\mathcal{C}_{i+1}$ is as required.

  Suppose that $C_{v_{i}}$ contains two vertices $v_i$ and $v_i'$ which have degree 3 in $\mathcal{C}_i'$.
 Since $v_i' \in V(C_{v_i})$ and $v_i'\not=v_i$, if $v_i'$ is a dominated vertex, then by Claim \ref{cvi}, $v_i'$ is dominated by a vertex in $C_{v_i}$. 
 Recall the subdivision operation on $G_i/C_{v_i}$ as shown in Figure \ref{cs}: 
\begin{enumerate}[(i)]
\setlength{\itemsep}{1pt}
\setlength{\parsep}{1pt}
\setlength{\parskip}{1pt}
     \item 
 If $v_i$ ($v_i'$) is   dominated by some vertex in $V(C_{v_i})$, then the edge $e \in E(\mathcal{C}_i')\backslash E(C_{v_i} )$ incident with $v_i$ ($v_i'$) is subdivided once, and the subdivided vertex is dominated by $v_i$ in $\mathcal{C}_{i+1}$;
 \item If
 $v_i$  is  dominated by its predecessor $y_i$ in $P_{u_i}$, then by ($*$), 
 the bad path $y_iz_{i1}z_{i2}v_i \ss \mathcal{C}_{i+1}$, and $z_{i1},z_{i2}$ are contracted as $\mathcal{C}_{i+1}$ becomes $\mathcal{C}_i'$;

 \item
 If $v_i$ ($v_i'$) is a dominating vertex, then the edge $e \in E(\mathcal{C}_i')\backslash E(C_{v_i} )$ incident with $v_i$ ($v_i'$)  is not subdivided as $G_i$ becomes $G_{i+1}$. 
\end{enumerate}

If none of $\{v_i,v_i'\}$ are dominating vertices in $G_i$, then by (i), (ii) and Claim \ref{cvi}, we can easily check that $\mathcal{C}_i=\mathcal{C}_i'$ is required. Thus, we are left to consider the case when at least  one of $\{v_i,v_i'\}$ is a dominating vertex in $G_i$. 
Let $P',P''$ be two internal disjoint $(v_i,v_i')$-path in $C_{v_i}$. In this case, we must remove all the internal vertices of $P'$ or $P''$ from $\mathcal{C}_i'$ to get $\mathcal{C}_i$ as expected. We now choose a path $P'$ to be kept in $\mathcal{C}_{i}'$ as follows.
If there exists a bad $(y_j,v_j)$-path $y_jz_{j1}z_{j2}v_j \ss C_{v_i}$, then $v_j\not=v_i$ by Claim \ref{detail} since $x$ lies in $\mathcal{C}_{i+1}$. If $v_j=v_i'$,  we assume without loss of generality that $y_jz_{j1}z_{j2}v_j\ss P'$. In this case, if $v_i$ is dominated by a vertex in $C_{v_i}$, then by the construction of $P_{u_i}$ and Claim \ref{detail}, $v_i$ is dominated by its successor in $P'$. 
If one of $\{v_i,v_i'\}$ is dominated by a vertex $w$ in $C_{v_i}$, then note that the other one is a dominating vertex, we assume $P'$ to be the path containing $w$. 
For other cases, we choose one of  $(v_i,v_i')$-path in $C_{v_i}$ as $P'$ arbitrarily.
Let $\mathcal{C}_i=\mathcal{C}_i'- (V(C_{v_{i}})\backslash V(P'))$. By the choice of $P'$ and (i)-(iii), we can check that $\mathcal{C}_i$ is an alternative subgraph containing $x$ and at least another dominating vertex in $\{v_i,v_i'\}$,
and satisfies ($*$).

By the arguments above, we can obtain an alternative subgraph $\mathcal{C}_i$ containing $x$ and at least two dominating vertices in $G_i$ for all $1 \leq i \leq t$, and $\mathcal{C}_1=H_0$ is an alternative subgraph of $G$ as required. \bqed

\vskip 2mm
By the graph as shown in Figure \ref{tight}, the lower bound is tight for the number of dominating vertices contained in the alternative  subgraph $H_0$ in Lemma \ref{vrh2}.

\vskip 2mm
At the end of this section, we give the following lemma which acts as the inductive basis when showing our main results by induction on $|D|$.

\begin{lem}\label{r=1}
    If $|D|=1$ in $(G,D)$, then there exists an orientation $\rbjt{G}$ of $G$ such that $diam(\rbjt{G})\leq 4$ and $sdiam(\rbjt{G})\leq 6$.
\end{lem}

\noindent{\bf{Proof}:} Let $v$ be the unique dominating vertex.
We give the orientation $\rbjt{G}$ of $G$ as follows.  We first orient as many as possible  edge-disjoint triangles containing $v$ as directed triangles. Then repeat the following orientations for a dominated vertex $v'$ not in a directed triangle containing $v$: Since $G$ is bridgeless and $r=1$,  there must exist a dominated vertex $v''$ adjacent to $v'$, and $v''v$ has been oriented.
Orient $vv'v''v$ as a directed triangle according to the orientation of $v''v$. By the above orientations, we can see that each dominated vertex lies in a directed triangle containing $v$. Hence,  $diam(\rbjt{G})\leq 4$ and $sdiam(\rbjt{G})\leq 6$. \bqed

\section{Oriented diameter and domination number}
\label{3}

\noindent{\bf{\emph{Proof of Theorem  \ref{main}}}:}  Let $D$ be a dominating set of $G$ with $|D|=r$, and $(G,D)$ be a standard pair.  We use induction on $r$. By Lemma \ref{r=1}, Theorem \ref{main} is true for $r=1$.

Suppose that $r \geq 2$ and Theorem \ref{main} holds for small $r$. Before proceeding on, we first establish a claim to facilitate calculating distances. 
For a bridgeless subgraph $H\ss G$ with a strong orientation $\rbjt{H}$, let 
$r'=|D\cap V(H)|$, and
$d_i(\rbjt{H})$ denote the  maximum directed distance between any two vertices $u,v \in V(H)$ with $|\{u,v\}\cap D|=i$. Set $$m(\rbjt{H})=max\{d_0(\rbjt{H})-2, d_1(\rbjt{H})-1, d_2(\rbjt{H})\}.$$ 

\begin{claim}\label{contract}  Let  $H\ss G$ be a bridgeless subgraph in which every dominated vertex is dominated by a vertex in $ V(H)$ and $r'\geq 2$. Then for any strong orientation $\rbjt{H}$, $\rbjt{H}$ can be extended to a strong orientation $\rbjt{G}$ such that  
 $diam(\rbjt{G}) \leq \lceil \frac{7(r-r'+1)+1}{2}\rceil+m(\rbjt{H})$.
\end{claim}

\noindent{\bf{Proof.}} 
Let $h$ be the contracted vertex in $G/H$, $G_1$ be $S(G/H)$ with a new isolated triangle associated with $h$   and $D_1=D\cup \{h\}\backslash V(H)$. Clearly, $(G_1,D_1)$ is a standard pair with $|D_1|=r-r'+1<r$.
By induction hypothesis, $G_1$ admits an orientation $\rbjt{G}_1$ such that $diam(\rbjt{G}_1) \leq \lceil \frac{7(r-r'+1)+1}{2}\rceil$.
Note that $G/H$ becomes $G_1$ after the subdivision operation. Denote $\rbjt{G/H}$ as the orientation of $G/H$ corresponding to $\rbjt{G}_1$. Combining  $\rbjt{G/H}$ and $\rbjt{H}$, we can obtain an orientation $\rbjt{G}$ of $G$.
  
Now we calculate the directed distance  between any two vertices $u,v$ in $\rbjt{G}$.

If $u,v \in V(H)$, then 
$$\theta(u,v) \leq diam(\rbjt{H})\leq 2+m(\rbjt{H}) < \left\lceil \frac{7(r-r'+1)+1}{2}\right\rceil+m(\rbjt{H}).$$

If $u\notin V(H)$ and $v \in V(H)$, then note that $\theta_{\rbjt{G/H}}(u,h) \leq \theta_{\rbjt{G}_1}(u,h)$, and $d(\rbjt{G}_1)=2$ if $G_1$ is a triangle, by Lemma \ref{-2} we have
 \begin{align*}
\theta(u,v)\leq \theta_{\rbjt{G/H}}(u,h)+diam(\rbjt{H})  
& \leq \theta_{\rbjt{G}_1}(u,h)+diam(\rbjt{H}) \\
& \leq
\left \lceil \frac{7(r-r'+1)+1}{2}\right\rceil-2+diam(\rbjt{H}) \\
& \leq \left\lceil \frac{7(r-r'+1)+1}{2}\right\rceil +m(\rbjt{H}).
 \end{align*}

Suppose that $u, v \notin V(H)$. If there exists a shortest directed $(u,v)$-path   not via $h$ in $\rbjt{G/H}$, then $$\partial(u,v) = \partial_{\rbjt{G/H}}(u,v) \leq \partial_{\rbjt{G}_1}(u,v) \leq 
\left\lceil \frac{7(r-r'+1)+1}{2}\right\rceil<\left\lceil \frac{7(r-r'+1)+1}{2}\right\rceil+m(\rbjt{H}).$$

If any shortest directed path from $u$ to $v$ must via $h$ in $\rbjt{G/H}$, denote one of such directed path as $\rbjt{P}=u\cdots h \cdots v$. Combining with $\rbjt{H}$, there is a directed ($u,v$)-path $\rbjt{P_G}=u\cdots h_1\cdots h_2 \cdots v$
 in $\rbjt{G}$, where $\rbjt{P_G} \cap \rbjt{H}=h_1\cdots h_2$.
 By the subdivision operation,  if  both $h_1,h_2$ are dominated vertices, then 
 $$\partial(u,v) \leq \partial_{\rbjt{G/H}}(u,v)+d_0(\rbjt{H}) \leq \partial_{\rbjt{G}_1}(u,v)-2+d_0(\rbjt{H}) \leq \left\lceil \frac{7(r-r'+1)+1}{2}\right\rceil+m(\rbjt{H});$$
If  one of $\{h_1,h_2\}$ is a dominating vertex and another is a dominated vertex, then 
$$\partial(u,v) \leq \partial_{\rbjt{G/H}}(u,v)+d_1(\rbjt{H}) \leq \partial_{\rbjt{G}_1}(u,v)-1+d_1(\rbjt{H}) \leq \left\lceil \frac{7(r-r'+1)+1}{2}\right\rceil+m(\rbjt{H}); $$
If  both $h_1,h_2$ are dominating vertices, then 
$$\partial(u,v) \leq \partial_{\rbjt{G/H}}(u,v)+d_2(\rbjt{H}) \leq \partial_{\rbjt{G}_1}(u,v)+d_2(\rbjt{H}) \leq \left\lceil \frac{7(r-r'+1)+1}{2}\right\rceil+m(\rbjt{H}).$$ 
Hence, $diam(\rbjt{G}) \leq \lceil \frac{7(r-r'+1)+1}{2}\rceil+m(\rbjt{H})$.
\hfill{$\square$}

\vskip 5mm
Now, we continue our proof. By Lemma \ref{vrh2},  let $H_0$ be an alternative subgraph of $G$ containing maximum number of dominating vertices. And subject to these, $H_0$ is edge-minimal. Set $r_0=|V(H_0)\cap D|$. Then $r_0\ge 2$ and $|V(H_0)|=3r_0$.
Obviously, $r_0\leq r$. Assume that $\rbjt{H}_0$ is a strong orientation of $H_0$. Clearly, $m(\rbjt{H}_0)\leq diam(\rbjt{H}_0)\le |H_0|-1=3r_0-1$.

 If  $r_0 \geq 5$, by applying Claim \ref{contract} on $\rbjt{H}_0$, there is an orientation $\rbjt{G}$ of $G$ such that 
$$diam(\rbjt{G}) \leq \left\lceil \frac{7(r-r_0+1)+1}{2}\right\rceil+3r_0-1=\left\lceil \frac{7r+1-r_0+5}{2}\right\rceil \leq \left\lceil \frac{7r+1}{2}\right\rceil,$$ 
and hence Theorem \ref{main} holds. 
\vskip 2mm

If $r_0=3,4$, then $|H_0|=9$ or 12. By analyzing the structures of $H_0$ based on the definition of alternative subgraph $H_0$ and the minimality of $H_0$, we can find a strong orientation $\rbjt{H}_0$ such that $m(\rbjt{H}_0)$ is less than $|H_0|-1$.

\begin{claim}\label{r_034}
    If $ r_0=3, 4$, then $H_0$ admits a strong orientation $\rbjt{H}_0$ with $m(\rbjt{H}_0)\leq 3r_0-2$.
\end{claim}
\noindent{\bf{Proof.}}
 Recall that $m(\rbjt{H}_0)=max\{d_0(\rbjt{H}_0)-2, d_1(\rbjt{H}_0)-1, d_2(\rbjt{H}_0)\}$, and for any strong orientation $\rbjt{H}_0$ of $H_0$, $\max\{d_0(\rbjt{H}_0)-2, d_1(\rbjt{H}_0)-1\} \leq 3r_0-2$. We remain to find a strong orientation $\rbjt{H}_0$ of $H_0$ such that $d_2(\rbjt{H}_0) \leq 3r_0-2$. If  there exists no Hamiltonian path in $H_0$ between any two dominating vertices, then $d_2(\rbjt{H}_0)\leq 3r_0-2 $ for any strong orientation $\rbjt{H}_0$. So we suppose that there exists a Hamiltonian path $P=u_1\cdots u_{3r_0}$ in $H_0$, where $u_1,u_{3r_0}$  are dominating vertices, and $u_2,u_{3r_0-1}$ is dominated by $u_1,u_{3r_0}$, respectively. 
 Since $H_0$ is an alternative subgraph, denote $u_i,u_j$ as another  vertex dominated by $u_1,u_{3r_0}$, respectively, $3\leq i,j \leq 3r_0-2$. 
 
 For an edge $u_su_t\in E(H_0)$ with $s<t$, we say an orientation of $u_su_t$ is \emph{forward} if  $u_su_t$ is oriented from $u_s$ to $u_t$.
For the cycle $C=u_sPu_tu_s \ss H_0$ and $s<t$, we say a strong orientation of $C$ is  \emph{forward}  if   $u_sPu_t$ is oriented as a directed path from $u_s$ to $u_t$. Otherwise, we say the orientation of $u_su_t$ or $C$ is  \emph{backward}.

Now, we prove that $H_0$ admits a strong orientation $\rbjt{H}_0$ such that the shortest directed path between any two dominating vertices is not Hamiltonian as follows.

If $i>j$, then by the minimality of $H_0$, $H_0$ is $P$ plus two edges $u_1u_i, u_ju_{3r_0}$. Since $H_0$ is alternative, we can check that $(i,j)=(7,3)$ if $r_0=3$, and $(i,j)=(10,3),(10,6)$ or $(7,3)$ if $r_0=4$. Because the two graphs corresponding to $(i,j)=(10,6)$ and $(7,3)$ are isomorphic, we assume $(i,j)=(10,3),(7,3)$ if $r_0=4$. Orient $u_1Pu_iu_1$  and $u_ju_{3r_0}$ forward, and all the edges of $u_iPu_{3r_0}$ backward, we get the required $\rbjt{H}_0$. Therefore, we are left to consider the case when $i<j$.

Assume that  $H_0$ has  at least two cut vertices.
If $(i,j)=(3,3r_0-2)$, then by the minimality of $H_0$, $u_3$ and $u_{3r_0-2}$ are two cut vertices, and the two triangles $u_1u_2u_3u_1$ and $u_{3r_0-2}u_{3r_0-1}u_{3r_0}u_{3r_0-2}$ are blocks of $H_0$. 
If $(i,j)\not=(3,3r_0-2)$, then since $H_0$ is alternative, we can deduce $r_0=4$, and $(i,j)=(3,7), (6,10)$. By the symmetry of $u_1$ and $u_{12}$, assume $(i,j)=(3,7)$. In this case, $u_3$ is a cut vertex. Note that the possible cut vertex must be a dominated vertex among $u_3,u_4,u_6,u_7$, $H_0$ has exactly two cut vertices and the other one is $u_6$ or $u_7$. If $u_7$ is a cut vertex, then  by the minimality of $H_0$,
$u_1u_2u_3u_1$ and $u_7Pu_{12}u_7$ are two blocks. If $u_6$ is a cut vertex, then we can check $u_3u_6\in E(H_0)$. Hence, $u_1u_2u_3u_1$ and $u_3Pu_6u_3$ are two blocks. By the above arguments, $H_0$ always has two blocks which are cycles. Orient one of the two cycles forward, another cycle backward, and give any other block an arbitrary strong orientation, we get the required  $\rbjt{H}_0$.

If  $u_iu_j \in E(P)$, then since $H_0$ is alternative, by symmetry of $u_1$ and $u_{r_0}$, we may assume that $(i,j)=(3,4)$  when $r_0\in \{3,4\}$, and $(i,j)=(6,7)$ when $r_0=4$. Since $H_0$ is bridgeless,  $H_0-u_iu_j$ has an edge $e$ between two vertex sets $\{u_1,\ldots,u_i\}$ and $\{u_j,\ldots,u_{3r_0}\}$.
By the minimality of $H_0$,  $H_0$ is $P$ plus three edges $u_1u_i,u_ju_{3r_0}$ and $e$. 
Note that if $r_0=4$ and $(i,j)=(6,7)$, then one endpoint of $e$ is contained in $\{u_2,u_3,u_5,u_6\}$. By the minimality of $H_0$, $u_5$ can not be one endpoint of $e$.
If $r_0=4$,  $(i,j)=(6,7)$ and $e\in \{u_2u_7,u_2u_{10}\}$, we  orient  $u_1Pu_iu_1$, $u_jPu_{12}u_j$ and $e$ forward,  and $u_iu_j$ backward. Otherwise, no matter $r_0=3$ or 4, we always orient the cycle $u_1Pu_iu_1$ and  $u_iu_j$ forward, and $u_jPu_{3r_0}u_j$ and $e$ backward. Then we get the required  $\rbjt{H}_0$.

If  $H_0$ has at most one cut vertex and $u_iu_j \notin E(P)$, then we can check $r_0=4$, and by symmetry we may assume $(i,j)=(3,7)$. Note that $u_3$ must be a cut vertex in $H_0$, there exists an edge $u_3u_s$. By the minimality of $H_0$, $s \in \{6,10,11\}$. 
If $s=6$, since $H_0$ is bridgeless, $H_0-u_6u_7$ has an edge $u_tu_{t'}$ with $t\in \{3,4,6\}$ and $t'\in \{7,8,10,11\}$. Since  $H_0$ has at most one cut vertex and is edge-minimal,  $H_0$ must be $P$ plus four edges  $u_1u_3, u_7u_{12},u_3u_6,u_tu_{t'}$ with $t=4$ and $t' \in \{10,11\}$.
In this case, orient $u_1u_2u_3u_1$, $u_3Pu_6u_3$, $u_7Pu_{12}u_7$ and $u_4u_{t'}$ forward, and  $u_6u_7$ backward.
If $s\in \{10,11\}$, then by the minimality of $H_0$, $H_0$ must be $P$ plus three edges $u_1u_3, u_7u_{12}, u_3u_s$. In this case, orient $u_1u_2u_3u_1$, $u_7Pu_{12}u_7$ and $u_3u_s$ forward, and all edges of $u_3Pu_7$ backward. 
By the above orientations, we get the required  $\rbjt{H}_0$.
\hfill{$\square$}

\vskip 2mm
If $r_0=3,4$, then by Claim \ref{r_034}, we may assume that $\rbjt{H}_0$ is a strong orientation such that $m(\rbjt{H}_0)\leq 3r_0-2$.   Applying Claim \ref{contract} on this $\rbjt{H}_0$,  there is an orientation $\rbjt{G}$ of $G$ such that $diam(\rbjt{G}) \leq \lceil \frac{7(r-r_0+1)+1}{2}\rceil+m(\rbjt{H}_0)\leq \lceil \frac{7r+1}{2}\rceil$, and hence Theorem \ref{main} holds.

\vskip 2mm
Finally, we consider the case $r_0=2$.  
By Lemma \ref{vrh2}, every $v\in D$ is contained in an alternative subgraph which contains exactly two dominating vertices. After an easy analysis, such alternative subgraph is isomorphic to    $H_{01}$ or $H_{02}$, as shown in Figure \ref{h2}. 

\begin{figure}[!htb]
\centering
\includegraphics[height=0.2\textwidth]{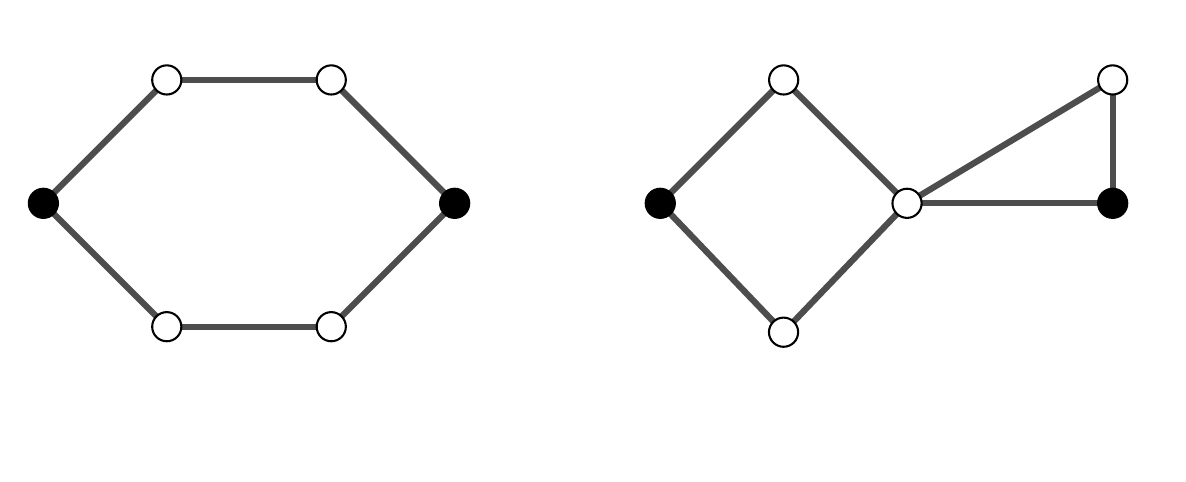}
\put(-162,10){\makebox(3,3){$H_{01}$}}
\put(-50,10){\makebox(3,3){$H_{02}$}}
\caption{Two $H_0$  with $r_0=2$ and minimal edges.}
\label{h2}
\end{figure}

Let $H_1,H_2$ be any two distinct subgraphs isomorphic to $H_{01}$ or $H_{02}$. 
Observe that  $|V(H_1)\cap V(H_2)\cap D|\leq  2$.
If $|V(H_1)\cap V(H_2)\cap D|= 1$, then since $r_0=2$, $H_1$ and $H_2$ have at most one common edge. Moreover, if $H_1$ and $H_2$  have exactly one common edge, then it is not difficult to check the union $H$ of $H_1$ and $H_2$ must be isomorphic to one of the six graphs shown in Figure \ref{h21}, and $m(\rbjt{H}) \leq 7$ for the orientations $\rbjt{H}$ shown in Figure \ref{h21}.
Applying Claim \ref{contract} on $\rbjt{H}$,  there is an orientation $\rbjt{G}$ of $G$ such that $diam(\rbjt{G}) \leq \lceil \frac{7(r-3+1)+1}{2}\rceil+m(\rbjt{H})\leq \lceil \frac{7r+1}{2}\rceil$. 

\begin{figure}[!htb]
\centering
\includegraphics[height=0.35\textwidth]{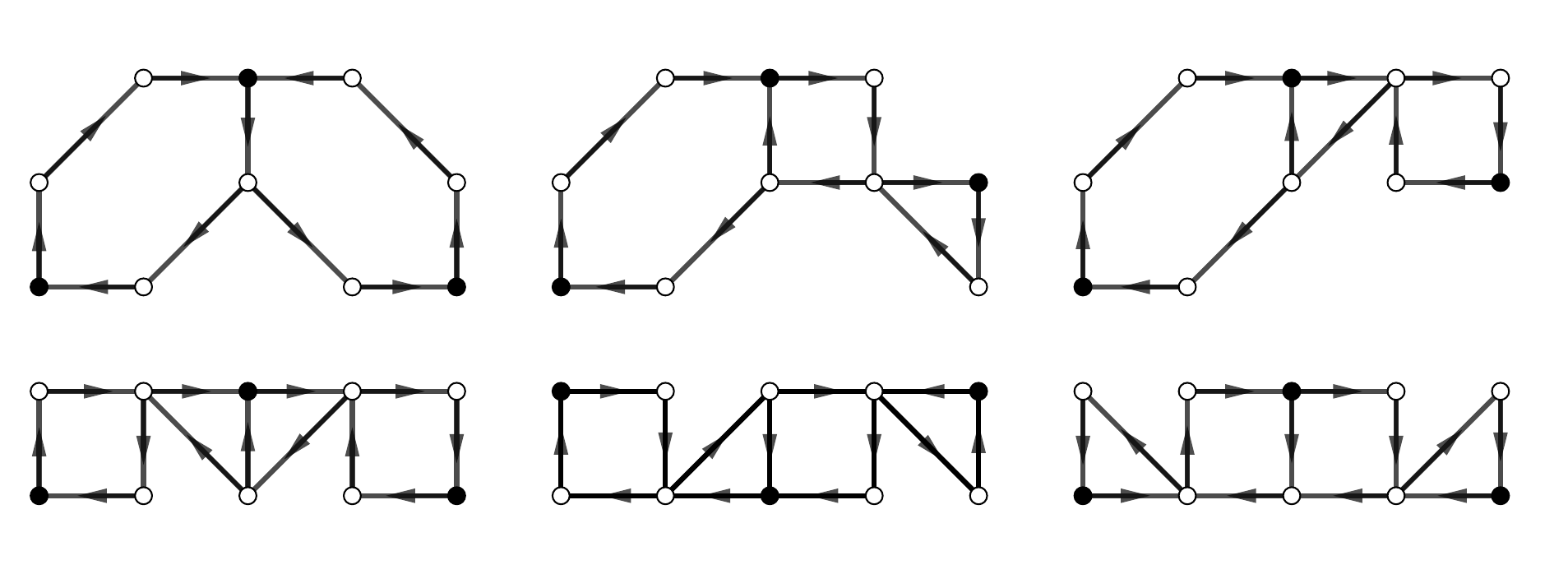}
\caption{Subgraphs $H$  with their orientations.}
\label{h21}
\end{figure}

Now, we may assume that if $|V(H_1)\cap V(H_2)\cap D|= 1$, then   $H_1,H_2$ are edge-disjoint, where $H_1,H_2$ are two distinct subgraphs isomorphic to $H_{01}$ or $H_{02}$.
Note that $H_1,H_2$ are also edge-disjoint if $|V(H_1)\cap V(H_2)\cap D|= 0$.
With these in mind,  
we construct a new graph $\mathcal{G}$ as follows:  $V(\mathcal{G})=D$, and  for any  $u,v\in V(\mathcal{G})$, $uv \in E(\mathcal{G})$  if and only if  there exists at least one subgraph  of $G$ isomorphic to $H_{01}$ or $H_{02}$ and containing both $u$ and $v$.
Obviously, $\mathcal{G}$ is a forest without isolated vertex.
For convenience, we say $uv \in E(\mathcal{G})$  corresponding    to $H_{01}$ or $H_{02}$. By the construction, two distinct edges of $\mathcal{G}$ correspond to two edge-disjoint subgraphs isomorphic to $H_{01}$ or $H_{02}$ in $G$, no matter the two edges are adjacent or not. We also say  a subgraph $\mathcal{H}\ss \mathcal{G}$ corresponding to a subgraph $H\ss G$.

If there exists a vertex of degree at least 3 in $\mathcal{G}$,  then there is a  star $\mathcal{S}_3\ss \mathcal{G}$ with a center of degree three. Let $H\ss G$ corresponding to  $S_3$, and $\rbjt{H}$ its  arbitrary strong orientation. Obviously, $m(\rbjt{H}) \leq diam(\rbjt{H})\leq 10$.
Applying Claim \ref{contract} on  $\rbjt{H}$,  there is an orientation $\rbjt{G}$ of $G$ such that $diam(\rbjt{G}) \leq \lceil \frac{7(r-4+1)+1}{2}\rceil+m(\rbjt{H})\leq \lceil \frac{7r+1}{2}\rceil-10+10 =\lceil \frac{7r+1}{2}\rceil$. 
Thus, we may assume the maximum degree of $\mathcal{G}$ is at most two, that is,
$\mathcal{G}$ is a linear forest consisting of vertex disjoint paths.

Let $\mathcal{P}$ be a longest path in $\mathcal{G}$, $H\ss G$ corresponding to  $\mathcal{P}$,  $\rbjt{H}$ a strong orientation of $H$, and $r'=|D\cap V(H)|\geq 2$. 
We need the following lemma to go on.

\begin{lem}\label{420}
(Kurz and  L\"atsch \cite{KL})
    If $H$ is a bridgeless subgraph of $G$ containing $D$, and $\rbjt{H}$ is its strong orientation, then there exists an orientation $\rbjt{G}$ of $G$ such that $diam(\rbjt{G})\leq \max\{d_0(\rbjt{H}),d_1(\rbjt{H})+2,d_2(\rbjt{H})+4\}$.
\end{lem}

\vskip 1mm

If $D\ss V(H)$, then $r=r'$. By the structure of $H$, we can calculate that
$d_0(\rbjt{H})\leq  5+\lceil \frac{7(r-3)}{2}\rceil +5 <  \lceil \frac{7r+1}{2}\rceil,$
$d_1(\rbjt{H})\leq  4+\lceil \frac{7(r-3)}{2}\rceil +5= \lceil \frac{7r+1}{2}\rceil-2,$ and
$d_2(\rbjt{H}) \leq  \lceil \frac{7(r-1)}{2}\rceil =  \lceil \frac{7r+1}{2}\rceil-4$. Thus, by Lemma \ref{420}, 
there is an orientation $\rbjt{G}$ of $G$ such that $diam(\rbjt{G})\leq \max\{d_0(\rbjt{H}),d_1(\rbjt{H})+2,d_2(\rbjt{H})+4\} \leq \lceil \frac{7r+1}{2}\rceil$.
\vskip 1mm
In the rest of the proof, we will show that the case $D\backslash V(H)\not= \emptyset$ can not occur.
Otherwise,
    let $h$ be  the contracted vertex in $G/H$,  
    $G_1$ be $S(G/H)$ with a new isolated triangle associated with $h$,     
    and $D_1=D\cup \{h\}\backslash V(H)$. Clearly, $(G_1,D_1)$ is a standard pair with $|D_1|=r-r'+1<r$. Since  $D\backslash V(H)\not= \emptyset$, $|D_1|\geq 2$. Applying Lemma \ref{vrh2} on $(G_1,D_1)$ and $h$,
there exists an alternative subgraph $H_1\ss G_1$ containing $h$ and at least two dominating vertices. 
Let $H'\ss G/H$ and $H''\ss G$ be the subgraph corresponding to $H_1$, respectively.  
Since $d_{H_1}(h)=2$, by the subdivision operation,  we have $d_{H'}(h)=2$. Let $N_{H'}(h)=\{t_1,t_2\}$,
and $h_1t_1,h_2t_2\in E(H'')$ be two edges corresponding to $ht_1,ht_2\in E(H')$, respectively. Note that $H\ss H''$.

If $h_1=h_2$ and $h_1\in V(H)\cap D$, then 
by the subdivision operation,  $H_1$ is also an alternative subgraph of $G$ containing at least two dominating vertices. Since $r_0=2$, $H_1$ is isomorphic to $H_{01}$ or $H_{02}$. Thus, $H''$ is the union of $H$ and $H_1$ with one dominating vertex $h_1$ in common. But this contradicts that $\mathcal{G}$
is a linear forest, or
the maximality of $\mathcal{P}$.

If $h_1=h_2$ and $h_1\in V(H)\backslash D$, we write the $H_{01}$ or $H_{02}$ in $H$ containing $h_1$ as $H_0$.
Then by
the subdivision operation, we can check that the union of this $H_0$ and $H'$ with one dominated vertex $h_1$ in common, is an alternative subgraph of $G$ containing at least 3 dominating vertices, which contradicts $r_0=2$. 

Assume that $h_1\not=h_2$. 
Recall that $\odot$ and $\circledast$ denote a dominated and  a  dominating vertex  in $H$, respectively. If $h_1,h_2\in V(H)\cap D$,
by the structure of $H$, 
we can always find a path $P'\ss H$ of type $h_1\!\odot\!\odot\!\circledast\!\odot\!\cdots\! \odot\!\circledast\!\odot\!\cdots\!\odot\!\circledast\!\odot\!\odot h_2$. By the subdivision operation,  the union of $P'$ and $H'$ is an alternative subgraph of $G$ containing at least 3 dominating vertices, which contradicts $r_0=2$.

If  $h_1,h_2\in V(H)\backslash D$ and lie in the same $H_{01}$ or $H_{02}$, we write the $H_{01}$ or $H_{02}$ in $H$ containing $h_1,h_2$ as $H_0$.
Then by
the subdivision operation, we can check that the union of this $H_0$ and $H'$ with two dominated vertices $h_1,h_2$ in common, is an alternative subgraph of $G$ containing at least 3 dominating vertices, which contradicts $r_0=2$. 
If $h_1,h_2\in V(H)\backslash D$ and  lie in different $H_{01}$ or $H_{02}$, by the structure of $H$, we can always find a path $P'\ss H$ of type
$h_1\!\circledast\!\odot\cdots \odot\!\circledast\!\odot\cdots\odot\!\circledast h_2$, or a walk $P'\ss H$ in which each edge occurs once and of type $h_1\!\circledast\!\odot h_1 \!\odot\!\circledast\!\odot\cdots \odot\!\circledast\!\odot\cdots\odot\!\circledast h_2$ (only $h_1$ occurs twice and $h_1\!\circledast\!\odot h_1$ is a triangle in some $H_{02}$), or a walk $P'\ss H$ in which each edge occurs once and of type $h_1\!\circledast\!\odot h_1 \!\odot\!\circledast\!\odot\cdots \odot\!\circledast\!\odot\cdots \!\odot\!\circledast\!\odot \!h_2\! \odot\!\circledast h_2$ (only $h_1$ and $h_2$ occur twice and $h_1\!\circledast\!\odot h_1$, $h_2 \odot\!\circledast h_2$ are two triangles).
By the subdivision operation, the union of $P'$ and $H'$ is an alternative subgraph of $G$ containing at least 3 dominating vertices, which contradicts $r_0=2$.

Assume that $h_1\in V(H)\backslash D$ and $h_2\in V(H)\cap D$. If $h_1h_2\notin E(H)$, then by the structure of $H$, we can always find a path $P'\ss H$ of type 
$h_1\!\circledast\!\odot\!\cdots\! \odot\!\circledast\!\odot\!\cdots\!\odot\!\circledast\!\odot \!\odot h_2$, or a walk $P'\ss H$ in which each edge occurs once and of type $h_1\!\circledast\!\odot h_1\!\odot\!\circledast\!\odot \cdots \odot\!\circledast\!\odot\cdots\odot\!\circledast\!\odot \!\odot h_2$ (only $h_1$ occurs twice and $h_1\!\circledast\!\odot h_1$ is a triangle in some $H_{02}$). By the subdivision operation, the union of $P'$ and $H'$ is an alternative subgraph of $G$ containing at least 3 dominating vertices, which contradicts $r_0=2$.
If $h_1h_2\in E(H)$, then write the
$H_{01}$ or $H_{02}$ containing $h_1,h_2$ as $H_0$.
By the subdivision operation, if we replace the only subdivided vertex adjacent with $h$ in $H_1$ by $h_1$, the resulting graph $H_0'$ is an alternative subgraph of $G$ containing at least two dominating vertices, and shares exactly one edge $h_1h_2$ with $H_0$. Since $r_0=2$, $H_0'$ must be isomorphic to $H_{01}$ or $H_{02}$. But the existences of $H_0'$
and $H_0$ contradicts our assumption.
\bqed

\section{Oriented strong diameter and domination number}

{\noindent{\bf{\emph{Proof of Theorem \ref{main2}}}}:} Let $D$ be a dominating set of $G$ with $|D|=r$, and $(G,D)$ be a standard pair.  We use induction on $r$. By Lemma \ref{r=1}, Theorem \ref{main2} is true for $r=1$.

Assume that $r\ge 2$ and Theorem \ref{main2} holds for small values of $r$. By Lemma \ref{vrh2},  let $H_0$ be an alternative subgraph of $G$ containing maximum number of dominating vertices. Moreover, choose such $H_0$ to be edge-minimal. Set $r_0=|V(H_0)\cap D|$. Then $r_0\ge 2$ and $|V(H_0)|=3r_0$.
Obviously, $r_0\leq r$. 

Let $h_0$ be the contracted vertex in $G/H_0$, $G_1$ be $S(G/H_0)$ with a new isolated triangle associated with $h_0$, and $D_1=D\cup \{h_0\}\backslash V(H_0)$. Then $(G_1,D_1)$ is a standard pair with $|D_1|=r-r_0+1<r$. By induction hypothesis, $G_1$ admits an orientation $\rbjt{G}_1$ such that $sdiam(\rbjt{G}_1) \leq 7(r-r_0+1)-1$. Assume that $\rbjt{H}_0$ is a strong orientation of $H_0$. Note that $G/H_0$ becomes $G_1$ after the subdivision operation, we let $\rbjt{G/H_0}$ be the orientation of $G/H_0$ corresponding to $\rbjt{G}_1$. Since $sd(u,v)$ in $\rbjt{G/H_0}$ is no more than $sd(u,v)$ in $\rbjt{G}_1$ for any $u,v \in V(G/H_0)$, we have
\begin{align*}
\rbjt{sdiam}(G)&\le sdiam(\rbjt{G/H_0})+ sdiam(\rbjt{H}_0) \\
&\leq sdiam(\rbjt{G}_1)+ sdiam(\rbjt{H}_0) \leq 7(r-r_0+1)-1+sdiam(\rbjt{H}_0).
\end{align*}
If  $r_0 \geq 5$, then note that $sdiam(\rbjt{H}_0)\le 2(|H_0|-1)=2(3r_0-1)$, by the inequality above, we can deduce that $sdiam(\rbjt{G}) \leq 7(r-r_0+1)-1+2(3r_0-1) = 7r-r_0+4 \leq 7r-1,$ and so Theorem \ref{main2} holds. Thus, we are left to consider the case when $2\le r_0\le 4$.

Note that  $|V(H_0)|=3r_0$, each dominating vertex has degree 2, and  $H_0$ is edge-minimal, we can check $e(H_0) \leq 7(r_0-1)$ as follows. If $r_0=2$, it is easy to check that $H_0$ is isomorphic to one of the two graphs as shown in Figure \ref{h2}, and so $e(H_0)\le 7$. For $r_0=3$, if there exists one vertex of degree at least 4 in $V(H_0)\backslash D$, we have $e(H_0) \leq 14$ after a tedious check, and otherwise $e(H_0) \leq (3\times 6+2\times 3)/2<14$. So $e(H_0)\le 14$. For $r_0=4$, if there exists one vertex of degree at least 5 in $V(H_0)\backslash D$, then $e(H_0) \leq 21$ after a tedious check, and otherwise, $e(H_0) \leq (4\times 8+2\times 4)/2<21$, and hence $e(H_0)\le 21$. By the arguments above, $e(H_0) \leq 7(r_0-1)$. Thus, by the inequality above,  we have  $sdiam(\rbjt{G}) \leq 7(r-r_0+1)-1+e(H_0)\leq 7r-1$. \bqed

\section{Concluding Remarks}

One can see that the subgraph described in Lemma \ref{vrh2} plays a key role in proving Theorems \ref{main} and \ref{main2}, which provides us a recursive structure for induction. This structure may have potential applications for other problems, especially for the ones concerning domination number. 

Meanwhile, the upper bounds in Theorems \ref{main} and \ref{main2} can be attained by the graph as presented in Figure \ref{tight}, which we believe is the unique extremal graph subject to $G$ is edge-minimal and vertex-minimal. But the proof may be complicated. Note that the extremal graph has many cut vertices, it is of interesting to investigate the oriented diameter and oriented strong diameter of 2-connected graphs. 

\begin{pb}\label{p1} Is it true that $\rbjt{diam}(G)\leq 3\gamma-1$ and $\rbjt{sdiam}(G)\leq 3\gamma$ for
any 2-connected graph $G$ with $\gamma(G)= \gamma$? 
\end{pb}
If the answer to Problem \ref{p1} is affirmative, then the bounds are sharp as can be shown by a cycle on $3\gamma$ vertices.

\vskip 5mm
\noindent{\bf\large Acknowledgements}
\vskip 2mm

Wang was supported by NSFC under grant number 12401447, and Chen was supported by National Key R\&D Program of China under grant number 2024YFA1013900 and NSFC under grant number 12471327.

\end{document}